\input amstex

\documentstyle{amsppt}
\magnification1200
\pagewidth{6.5 true in}
\pageheight{9.25 true in}
\NoBlackBoxes

\topmatter
\title 
The distribution of smooth numbers in arithmetic progressions 
\endtitle
\author 
K. Soundararajan
\endauthor 
\address 
Department of Mathematics, 450 Serra Mall, Bldg. 380, Stanford University, 
Stanford, CA 94305-2125, USA
\endaddress
\email 
ksound{\@}stanford.edu
\endemail
\thanks The author is partially supported by the National Science Foundation (DMS 0500711) 
and the American Institute of Mathematics (AIM). 
\endthanks

\endtopmatter

\document

\head 1. Introduction \endhead

\noindent We say that a number $n$ is $y$-smooth if all the prime factors 
of $n$ lie below $y$.  Let ${\Cal S}(y)$ denote the 
set of all $y$-smooth numbers, and let ${\Cal S}(x,y)$ denote the 
set of $y$-smooth numbers below $x$.  Let $\Psi(x,y)$ denote the number 
of smooth integers below $x$; thus $\Psi(x,y)$ is the cardinality of 
${\Cal S}(x,y)$.  In this note we consider the distribution of 
smooth numbers among arithmetic progressions $a\pmod q$.  We suppose that 
$(a,q)=1$, and it is natural to expect that smooth numbers are equally 
distributed among such progressions: that is, 
$$
\Psi(x,y;q,a) := \sum\Sb n\in {\Cal S}(x,y) \\ n\equiv a\pmod q\endSb 
1 \sim \frac{1}{\phi(q)} \sum\Sb n\in {\Cal S}(x,y) \\ (n,q) =1 \endSb 1
=: \frac{1}{\phi(q)} \Psi_q(x,y). \tag{1}
$$
Naturally there are some limitations to when we may expect (1) to hold, but 
it seems safe to make the following conjecture. 

\proclaim{Conjecture I(A)}  Let $A$ be a given positive real number.  Let 
$y$ and $q$ be large with $q\le y^A$.  Then as $\log x/\log q \to \infty$ we 
have 
$$
\Psi(x,y;q,a) \sim \frac{1}{\phi(q)} \Psi_q(x,y).
$$
\endproclaim 

In [5], [6]  Granville established this Conjecture when $A<1$.  He noted that 
establishing the conjecture for arbitrarily large $A$ would be difficult, 
since that would imply Vinogradov's conjecture that the least 
quadratic non-residue $\pmod p$ lies below $p^{\epsilon}$.  For, if $p$ 
is a prime and $y$ lies below the least quadratic non-residue $\pmod p$ then 
all elements of ${\Cal S}(y)$ are quadratic residues $\pmod p$, and 
we cannot have the equidistribution property (1).  The best known 
result towards Vinogradov's conjecture is that the least quadratic 
non-residue $\pmod p$ lies below $p^{1/4\sqrt{e}}$.  Thus it would 
be interesting to establish Conjecture I(A) for $A<4\sqrt{e}$, and even more interesting 
to establish it for larger $A$.  In this 
context, Harman [9] has shown that for $q$ cube-free, $q\le y^{4\sqrt{e}-\epsilon}$ 
and $q^{2+\epsilon} \le x\le q^{1/\epsilon}$ one has $\Psi(x,y;q,a) 
\gg \Psi_q(x,y)/\phi(q)$.  A slightly weaker result holds for more 
general $q$; see also the work of Balog and Pomerance [2] in this 
direction.  

\proclaim{Theorem 1} Let $y$ and $q$ be large with $q\le y^{4\sqrt{e}-\epsilon}$.  
For $\exp(y^{1-\epsilon}) \ge x\ge y^{(\log\log y)^4}$ the asymptotic 
formula (1) holds. 
\endproclaim 

It seems plausible that with greater effort our methods could 
be extended to obtain Conjecture I when $A <4\sqrt{e}$.  We hope 
that an interested reader will accept that challenge.  

As remarked above, there is a serious obstacle to establishing 
Conjecture I for any larger value of $A$.  Namely, it may happen 
that the smooth numbers mostly lie in some subgroup of 
the group of reduced residues $\pmod q$; for example, the 
subgroup of quadratic residues.  However, within that subgroup we 
would expect equidistribution.  

\proclaim{Conjecture II(A)} Let $A$ be a given positive real 
number.  Let $y$ and $q$ be large with $q\le y^A$.  There 
exists a constant $C(A)$ depending only on $A$, and a subgroup 
$H$ of $({\Bbb Z}/q{\Bbb Z})^*$ of index at most $C(A)$ such that 
for any reduced residues $a$ and $b\pmod q$ with $a/b\in H$ 
we have, as $\log x/\log q \to \infty$ 
$$
\Psi(x,y;q,a) = \Psi(x,y;q,b) + o(\Psi_q(x,y)/\phi(q)). \tag{2}
$$
\endproclaim 

Conjecture I is the stronger statement that $H=({\Bbb Z}/q{\Bbb Z})^*$.  
We are optimistic that the methods developed here could be used 
to prove Conjecture II.  Towards that end, we prove the following 
Theorem.  

\proclaim{Theorem 2}  Let $A$ be any positive real 
number and $y$ and $q$ be large with $q\le y^A$.  There 
exists a subgroup $H$ of $({\Bbb Z}/q{\Bbb Z})^*$ of 
index bounded by $C(A)$ such that for any two residue 
classes $a$ and $b$ with $a/b\in H$ and all 
$\exp(y^{1-\epsilon}) \ge x \ge y^{(\log \log y)^4}$ 
the asymptotic formula (2) holds. 
\endproclaim

Let $a\pmod q$ be an arithmetic progression with $(a,q)=1$.  Using the 
orthogonality of the characters $\pmod q$ we may write 
$$
\Psi(x,y;q,a) = \frac{1}{\phi(q)} \sum_{\chi \pmod q} \overline{\chi(a)} \Psi(x,y;\chi),  
$$ 
where 
$$ 
\Psi(x,y; \chi) = \sum\Sb n\in {\Cal S}(x,y) \endSb \chi(n). \tag{3} 
$$ 
We expect that the main term arises from the principal character, 
and that the contribution of all other characters is negligible.   This is indeed 
the case for the range of Theorem 1.  In the range of Theorem 2 we shall 
establish that there are at most a bounded number (in terms of $A$) of characters (of bounded 
order)
for which the sum in (3) can be large.  The subgroup $H$ consists of those residue classes 
which take the value $1$ on all these problem characters.

{\bf Acknowledgments.}    The work on this paper was done while 
I was an Andr{\' e} Aisenstadt  Chair at the Centre de Recherches Math{\' e}matiques, 
Montr{\' e}al.   I am most grateful to CRM for their generous hospitality, and 
for providing a very stimulating environment.   I am also grateful to Andrew Granville 
for many valuable discussions on topics related to this paper. 

\head 2.  Preliminary Observations \endhead 

\noindent It is convenient to introduce a smooth weight $\Phi(x)$.   
We suppose that $\Phi:{\Bbb R}_{\ge 0} \to [0,1]$ is a function, smooth on that 
domain, and approximating the characteristic function of the interval $[0,1]$.   Concretely, 
given $\epsilon >0$ we shall take $\Phi$ to be $1$ on $[0,1-\epsilon]$, $0$ on $[1,\infty)$ 
so that $\Phi$ approximates from below the characteristic function of $[0,1]$ or 
we shall take $\Phi$ to be $1$ on $[0,1]$ and $0$ on $[1+\epsilon,\infty)$ 
getting an approximation from above.  At the last step, we shall let $\epsilon$ go to zero.    
With such a choice for $\Phi$, we shall consider 
$$
\Psi(x,y;q,a, \Phi) = \sum\Sb n \in {\Cal S}(y) \\ n\equiv q\pmod q \endSb  \Phi(n/x) = 
\frac{1}{\phi(q)} \sum_{\chi \pmod q} \overline{\chi(a)} \Psi(x,y;\chi,\Phi)  \tag{2.1}
$$
with 
$$
\Psi(x,y;\chi,\Phi) = \sum\Sb n\in {\Cal S}(y) \endSb \chi(n) \Phi(n/x). \tag{2.2}
$$

We define for Re$(s)>0$
$$
L(s,\chi;y) = \prod_{p\le y} \Big(1-\frac{\chi(p)}{p^s} \Big)^{-1} = \sum_{n\in {\Cal S}(y)} \frac{\chi(n)}{n^s}.
$$
By Mellin inversion we note that, for any $c>0$,  
$$
\Psi(x,y;\chi,\Phi) = \frac{1}{2\pi i} \int_{c-i\infty}^{c+i\infty} 
L(s,\chi;y) x^s {\check \Phi}(s) ds,  \tag{2.3} 
$$
where 
$$
{\check \Phi}(s) =\int_0^{\infty} \Phi(t) t^{s-1} dt.  \tag{2.4}
$$
Repeated integration by parts shows that for Re$(s)>0$, 
$$
{\check \Phi}(s) = -\int_0^{\infty} \Phi^{\prime}(t) \frac{t^s}{s} dt = \int_{0}^{\infty} 
\Phi^{\prime\prime}(t) \frac{t^{s+1}}{s(s+1)} dt = \ldots,
$$
so that, for any integer $k\ge 1$,
$$
|{\check \Phi}(s)| \ll_{\Phi,k} \frac{1}{|s|(|s|+1) \cdots (|s|+k-1)}. \tag{2.5}  
$$
In practice, we shall need (2.5) only for some fixed large number $k$; 
certainly $k=100$ will be sufficient. 


 Hildebrand and Tenenbaum [11] (see also the expository article [12]) developed the 
saddle point method to obtain an asymptotic for $\Psi(x,y)$.  Their results 
give, with some obvious modifications, an asymptotic formula for (2.3) in the 
case when $\chi =\chi_0$ is the principal character.  Let us begin by 
recalling some details of this result.  In order to keep our argument transparent, 
we will assume throughout that $\exp(y^{1-\epsilon}) \ge x\ge y^{(\log \log y)^4}$.  
Moreover, since Granville's work applies when $q\le \sqrt{y}$ we assume from now 
on that $\sqrt{y} \le q\le y^A$.   With more work we could relax these assumptions, and 
avoid the appeal to Granville's work.  

In the Hildebrand-Tenenbaum argument, the line of integration in (2.3) is 
chosen carefully.  They take $c$ to be $\alpha=\alpha(x,y)$ which is the unique solution to 
$$
\sum_{p \le y} \frac{\chi_0(p) \log p}{p^{\alpha}-1}  = \log x. 
\tag{2.6}
$$
As usual, we set $u=(\log x)/\log y$.  
For $y \ge (\log x)^{1+\epsilon}$ we have  
(see Lemmas 1 and 2 of [11])
$$
\alpha(x,y) = 1 -  \frac{\xi(u)}{\log y}+ O\Big(\frac{1}{\log x} +\frac{\log x}{y\log y}\Big), \tag{2.7a}
$$ 
where $\xi(u)$ is the unique solution to $e^{\xi} = 1+\xi u$ and it satisfies 
$$
\xi(u) \sim \log (u\log u). \tag{2.7b}
$$
Note that in our range for $x$ and $y$, we have that $\alpha \gg \epsilon$ is 
bounded away from zero.
With this choice for $c$, their asymptotic is 
$$
\Psi(x,y;\chi_0,\Phi) \sim \frac{x^{\alpha} L(\alpha,\chi_0;y){\check \Phi}(\alpha)}{ 
\sqrt{2\pi \phi_{2}(\alpha,\chi_0;y)}} \tag{2.8}
$$
where, in our range of $x$ and $y$,
$$
\phi_2(\alpha,\chi_0;y)= \sum_{p\le y} \chi_0(p) \frac{p^{\alpha}}{(p^{\alpha}-1)^2} 
\log^2 p \asymp \log x\log y. \tag{2.9}
$$ 
We also record that in our range for $x$ and $y$ we have 
$$
\log L(\alpha,\chi_0;y) \sim u; \tag{2.10} 
$$
this follows by a simple partial summation argument.

Take $c=\alpha$ in (2.3), and note that $|L(\alpha+it,\chi;y)| \le L(\alpha,\chi_0;y)$.  The 
rapid decay of ${\check \Phi}(s)$ (see (2.5)) allows us to truncate the integral in (2.3):
$$
\align
\Psi(x,y;\chi,\Phi) &= \frac{1}{2\pi i} \int_{\alpha - i\sqrt{q}}^{\alpha+ i\sqrt{q}} 
L(s,\chi;y) x^s {\check \Phi}(s) ds 
+ O( L(\alpha,\chi_0) x^{\alpha} q^{-10}) 
\\
&=  \frac{1}{2\pi i} \int_{\alpha - i\sqrt{q}}^{\alpha+ i\sqrt{q}} 
L(s,\chi;y) x^s {\check \Phi}(s) ds +O(\Psi(x,y;\chi_0,\Phi) q^{-2}).
\tag{2.11}
\\
\endalign
$$

To bound $\Psi(x,y;\chi,\Phi)$ for non-principal characters, we divide 
the characters $\pmod q$ into various sets based on the location of 
the zeros of $L(s,\chi)$.  Let $0\le j\le (\log q)/2$ be an 
integer, and let ${\Cal R}_j(q)$ denote the region $\{ s: \ \ 
\text{Re}(s) > 1- j/\log q, \ \ |\text{Im}(s)| \le q\}$.   The set $\Xi(j)$ 
consists of the non-principal characters $\chi$ for which 
$L(s,\chi)$ has no zeros in ${\Cal R}_j(q)$, but has a zero in
${\Cal R}_{j+1}(q)$.  By the log-free zero density estimate (see, for 
example, Chapter 18 of Iwaniec and Kowalski [13]) we know that 
$$
|\Xi(j)| \le C_1 e^{C_2 j}, \tag{2.12}
$$ 
for some absolute positive constants $C_1$ and $C_2$.  

There are three basic arguments used in the proof.  If $\chi \in \Xi(j)$ for 
some $j \ge 10A\log \log q$, a direct use of the implied zero-free region leads to a 
good bound for $\Psi(x,y;\chi,\Phi)$.  This takes care of the 
vast majority of characters.  
Second, in the region where $j \le 10A\log \log q$ but $j \ge 4A\log A+D$ for 
some absolute positive constant $D$,  we use a Rodosski\v i type argument (see 
[15]  and Chapter 9 of [14]) to bound $\Psi(x,y;\chi,\Phi)$.   
We are then left with a bounded number of problematic characters.  We 
show that those problem characters have bounded order, and the subgroup 
$H$ of Theorem 2 arises as the group of residue classes $r$ with $\chi(r) =1$ 
for all these problem characters.  Lastly when $A < 4\sqrt{e}$, Burgess's 
character sum estimates (see [3]) and  
reasoning along the lines of Vinogradov's $\sqrt{e}$ argument lead to the treatment
of problem characters, and thus to Theorem 1.

\head 3. Consequences of a zero-free region:  Basic argument \endhead

\proclaim{Lemma 3.1}  Let $\chi \pmod q$ be a non-principal character with 
$\chi \in \Xi(j)$ for some $j\ge 0$, and let $|t|\le q/2$.   Then for any $z\ge 2$ 
$$
\sum_{n\le z} \Lambda(n)\chi(n)n^{-it} \ll \frac{z(\log qz)^2}{q} + z^{1-j/\log q} (\log q)^2.
$$
\endproclaim 

\demo{Proof}  We may assume that $j\ge 1$ else the bound is 
trivial.   Therefore there are no issues with Siegel zeros.   
We follow a modification to the standard explicit formula argument (see 
for example Chapter 19 of Davenport [4]).  That argument shows 
$$
\sum_{n\le z} \Lambda(n) \chi(n) n^{-it} =-\sum_{|\gamma-t|\le q/2} \frac{x^{\rho-it}}{\rho-it} 
+O(z(\log qz)^2/q) + O(z^{\frac{1}{2}}), 
$$
where $\rho$ runs over the non-trivial zeros of $L(s,\tilde{\chi})$ 
with $\tilde{\chi}$ being the primitive character inducing $\chi$.  
Since $\chi \in \Xi(j)$ we see that if $\rho=\beta+i\gamma$ with $|t-\gamma|\le q/2$ 
then $|\gamma|\le q$ and so $\beta \ge 1-j/\log q$.  Using this, and splitting the 
sum over $\gamma$ into intervals of length $1$, and noting that each such 
interval has $\ll \log q$ zeros, we obtain the Lemma. 
\enddemo

\proclaim{Lemma 3.2}  Retain our assumptions on $x$, $y$ and $q$.  
Suppose that $j\ge 10A \log \log q$, and that 
$\chi$ is a non-principal character lying in $\Xi(j)$.   Let $\alpha=\alpha(x,y)$ be 
as in (2.6).  Let $B$ be a suitably large, but fixed positive number.  
For any $\sigma \ge \alpha -Bj/\log x -3 \log \log x/\log x$ and $|t|\le q/2$ we have 
$$
\log |L(\sigma+it,\chi;y)| =o(u).
$$
\endproclaim
\demo{Proof}  We split the primes below $y$ into the small primes $p\le q^{5(\log \log q)/j}$ 
and the large primes $q^{5(\log \log q)/j} < p \le y$.  Note that 
$q^{5(\log \log q)/j} \le q^{1/(2A)} \le \sqrt{y}$, by our assumptions on $j$, $y$ and $q$.
For the small primes we have 
$$
\align
\sum_{p\le q^{5(\log \log q)/j}} \log \Big| 1-\frac{\chi(p)}{p^{\sigma+it}}\Big|^{-1} 
&\ll \sum_{p\le q^{5(\log \log q)/j}} \frac{1}{p^{\alpha-Bj/\log x -3(\log \log x)/\log x}} \\
&\ll \sum_{p\le q^{5(\log \log q)/j}} \frac{1}{p^{\alpha}} \ll 
\sum_{p\le \sqrt{y}} \frac{1}{p^{\alpha}} \ll \sqrt{u}, \\
\endalign
$$
using the prime number theorem, partial summation, and (2.7a,b).

Next we treat the large primes.  Note that 
$$
\align
\sum_{q^{5(\log \log q)/j} \le p \le y}  \log \Big| 1-\frac{\chi(p)}{p^{\sigma+it}}\Big|^{-1} 
&=\sum_{q^{5(\log \log q)/j} \le n \le y} \text{Re } \frac{\chi(n)\Lambda(n)}{n^{\sigma+it}\log n} 
+ O\Big(\sum_{p\le y}\frac{1}{p^{2\sigma}} \Big).\\
\endalign
$$
The error term above is easily seen to be $o(u)$.  To handle the 
main term above we use partial summation together with Lemma 3.1.  Put temporarily 
$$
S(z) =\sum_{q^{5(\log \log q)/j }\le n\le z} \Lambda(n)\chi(n)n^{-it}.
$$
Partial 
summation gives
$$
\align
\sum_{q^{5(\log \log q)/j} \le n \le y} \text{Re } \frac{\chi(n)\Lambda(n)}{n^{\sigma+it}\log n} 
&=\text{Re }\int_{q^{5(\log \log q)/j}}^{y} \frac{1}{z^{\sigma} \log z} dS(z) 
\\
&\ll 1 +\frac{|S(y)|}{y^{\sigma} \log y} + \int_{q^{5(\log \log q)/j}}^{y}\frac{|S(z)|}{z^{\sigma+1}\log z}dz. 
\\
\endalign
$$
Using Lemma 3.1 we may check that $|S(z)|/z^{\sigma+1} \ll z^{-\alpha} (\log q)^{-3}$ 
in our range for $z$.  Using this above, along with (2.7a,b) we 
conclude that 
$$
\sum_{q^{5(\log \log q)/j} \le n \le y} \text{Re } \frac{\chi(n)\Lambda(n)}{n^{\sigma+it}\log n} 
\ll 
1+ \frac{u}{\log q} = o(u).
$$
We have established the desired bound for $\log |L(\sigma+it,\chi;y)|$. 
\enddemo

Suppose $\chi \in \Xi(j)$ for $j\ge 10 A \log \log q$.  With $B$ suitably large, we 
will apply Lemma 3.2 and shift the contour of integration in (2.11) to the 
$\alpha-Bj/\log x-3 \log \log x/\log x$ line.  By the 
rapid decay of ${\check \Phi}(s)$ the horizontal line segments contribute 
an amount $\ll \Psi(x,y;\chi_0,\Phi) q^{-2}$.   The remaining vertical 
line segment contributes, using Lemma 3.2, 
$$
\align
&\ll \int_{\alpha-Bj/\log x-3\log \log x/\log x-i\sqrt{q}}^{\alpha -Bj/\log x-3\log \log x/\log x+i\sqrt{q}} 
|L(s,\chi;y) x^s {\check \Phi}(s) ds| 
\ll e^{o(u)} x^{\alpha} e^{-Bj} (\log x)^{-3} 
\\
&\ll \Psi(x,y;\chi_0,\Phi) e^{-Bj} (\log x)^{-2},
\\
\endalign
$$
upon recalling (2.8), (2.9) and (2.10).   Thus, for $\chi\in \Xi(j)$ with 
$j\ge 10A \log \log q$ we conclude that 
$$
\Psi(x,y;\chi,\Phi) \ll \Psi(x,y;\chi_0,\Phi) \Big( \frac{e^{-Bj}}{(\log x)^2} + \frac{1}{q^2}\Big).  \tag{3.1}
$$

Now we can explain what suitably large means for $B$: namely that $B$ exceeds 
$C_2$, the constant appearing in the zero-density estimate (2.12).   Choosing $B$ that 
large, we conclude from (3.1) that 
$$
\sum_{j \ge 10A\log \log q} \sum_{\chi\in \Xi(j)} |\Psi(x,y;\chi,\Phi)| 
\ll \Psi(x,y;\chi_0,\Phi) \Big(\frac{1}{(\log x)^2} + \frac{1}{q}\Big). \tag{3.2}
$$
This is our basic zero-density argument, and it takes care of all but $\ll (\log q)^{10AC_2}$ 
characters $\chi \pmod q$.

\head 4.  Consequences of a zero-free region: The Rodosski\v i argument \endhead 

\noindent There remain $\ll (\log q)^{10AC_2}$ characters $\chi \pmod q$ which  are 
not covered by the argument of \S 3.  We now give a second argument to 
prune this set of characters, leaving only a bounded number of characters left to be estimated.   

\proclaim{Proposition 4.1}  Retain our ranges for $x$, $y$, and $q$.  
There exists an absolute positive constant $D$ such that 
if $\chi\in \Xi(j)$ with $j\ge 4A \log A+D$ then 
$$
\Psi(x,y;\chi,\Phi) \ll \Psi(x,y;\chi_0,\Phi)\left(  (\log x)e^{-\sqrt{u}/20} + q^{-2}\right).
$$
\endproclaim 

We shall bound $\Psi(x,y;\chi,\Phi)$ using (2.11).   Using (2.5),  we may 
express this bound as 
%
$$
\Psi(x,y;\chi,\Phi) \ll x^{\alpha}   \max_{|t|\le \sqrt{q}} |L(\alpha+it,\chi;y)|
+ q^{-2} \Psi(x,y;\chi_0,\Phi).  \tag{4.1} 
$$

We now define 
$$
{\Bbb D}_{\alpha}({ 1},\chi(p)p^{-it};y)^2 = \sum\Sb p\le y \\ p \nmid q \endSb 
\frac{1-\text{Re }\chi(p)p^{-it}}{p^{\alpha}}.
$$
This is a distance function which satisfies a triangle inequality: 
$$
{\Bbb D}_{\alpha}(f_1,g_1;y) + {\Bbb D}_{\alpha}(f_2,g_2;y) \ge 
{\Bbb D}_{\alpha}(f_1f_2, g_1g_2;y), 
$$
where $f_1$, $f_2$, $g_1$, $g_2$ are completely multiplicative functions 
taking values in the unit disc, and ${\Bbb D}_\alpha(f,g;y)^2 = \sum_{p\le y, p\nmid q} 
(1-\text{Re }\overline{f(p)}g(p))/p^{\alpha}$. 
The triangle inequality above may be deduced easily from Cauchy-Schwarz; see also 
the paper [8] for a general discussion of such inequalities, and [1], [7] for some 
applications.  
Note that 
$$
|L(\alpha+it,\chi;y)| \ll |L(\alpha,\chi_0;y)| 
\exp\Big( - {\Bbb D}_{\alpha}({ 1}, \chi(p)p^{-it};y)^2 \Big),
$$
and so from (4.1) (and recalling (2.8) and (2.9)) we obtain that 
$$
\align
\Psi(x,y;\chi,\Phi) &\ll x^{\alpha}   L(\alpha,\chi_0;y)  \exp\Big(-
\min_{|t|\le\sqrt{q}} {\Bbb D}_{\alpha}(1,\chi(p)p^{-it};y)^2 \Big)  
+  q^{-2} \Psi(x,y;\chi_0,\Phi) \\ 
&\ll \Psi(x,y;\chi_0,\Phi) \Big(\frac{1}{q^2} +( \log x) \exp\Big( -\min_{|t|\le\sqrt{q}} {\Bbb D}_{\alpha}(1,\chi(p)p^{-it};y)^2 \Big)  
\Big). \tag{4.2}\\
\endalign
$$


To proceed further we need some lower bounds on the distance function above; 
this is given in the following Lemma from which Proposition 4.1 is immediate.  
   
\proclaim{Lemma 4.2} Retain the notation of Proposition 4.1.   If 
$\chi \in \Xi(j)$ with $j\ge 4A \log A +D$ 
then, for $|t|\le \sqrt{q}$ we have
$$
{\Bbb D}_\alpha(1,\chi(p)p^{-it};y)^2 \ge \sqrt{u}/20. \tag{4.3}
$$
\endproclaim

The proof of Lemma 4.2 rests on some 
ideas of Rodosski\v i [15]; we follow here the treatment given in Chapter 9 of Montgomery [14].  
Observe 
that, using (2.7a,b),
$$
\align
{\Bbb D}_{\alpha}(1,\chi(p)p^{-it};y)^2 &\ge 
\sum\Sb \sqrt{y}\le p \le y\\ p\nmid q \endSb 
\frac{1-\text{Re }\chi(p)p^{-it}}{p^{\alpha}} 
\ge \frac{y^{\frac{1-\alpha}{2}}}{\log y} 
\sum\Sb \sqrt{y} \le p \le y\\ p\nmid q\endSb 
\frac{1-\text{Re }\chi(p)p^{-it}}{p} \log p\\
&\ge \frac{\sqrt{u}}{\log y} \sum\Sb \sqrt{y}\le p\le y \\ p\nmid q\endSb 
\frac{1-\text{Re }\chi(p)p^{-it}}{p}\log p. \tag{4.4} 
\\
\endalign
$$
Further let us define, as in Montgomery [14], the smooth weights 
$$
W(p) = 
\cases 
\log (p/\sqrt{y}) &\text{if  } \sqrt{y}\le p\le y^{\frac 34} \\ 
\log(y/{p}) &\text{if  } y^{\frac 34} \le p\le y\\ 
0&\text{otherwise}. \\
\endcases
$$
Note that for any $c>0$
$$
W(p) = \frac{1}{2\pi i} \int_{c-i\infty}^{c+i\infty} 
p^{-w} \Big(\frac{y^{w/2}-y^{w/4}}{w}\Big)^2 dw. 
$$
From (4.4), we see that 
$$
{\Bbb D}_{\alpha}(1,\chi(p)p^{-it};y)^2 \ge \frac{4\sqrt{u}}{\log^2 y} 
\sum\Sb p \nmid q \endSb \frac{1-\text{Re }\chi(p)p^{-it}}{p}  W(p)\log p. \tag{4.5} 
$$
 Since 
$$
\sum\Sb p \nmid q \endSb \frac{\log p}{p} W(p) 
\sim \frac{\log^2 y}{16},
$$
we may conclude from (4.5) that 
$$
{\Bbb D}(1,\chi(p)p^{-it};y)^2 \ge \frac{\sqrt{u}}{8} - \frac{4\sqrt{u}}{\log^2 y} 
\text{Re } \sum_{p} \frac{\chi(p)}{p^{1+it}} W(p) \log p. \tag{4.6} 
$$
The desired bound (4.3) is a consequence of the following Lemma.

\proclaim{Lemma 4.3} We keep the notations of Proposition 4.1. 
If $\chi \in \Xi(j)$ with $j\ge 4A \log A +D$ then, for $|t|\le \sqrt{q}$, we have 
$$
\text{Re } \sum_p \frac{\chi(p)}{p^{1+it}} W(p) \log p 
\le \frac{\log^2 y}{100}.
$$
\endproclaim 

\demo{Proof}  Let $\tilde{\chi}$ denote the primitive character inducing 
the character $\chi$.  Like $L(s,\chi)$, of course $L(s,\tilde{\chi})$ is also free of zeros 
in the region ${\Cal R}_j(q)$.
If $c>0$ then 
$$
\sum_p \frac{\chi(p) \log p}{p^{1+it}} W(p) = 
-\frac{1}{2\pi i} \int_{c-i\infty}^{c+i\infty} \frac{L^{\prime}}{L}(1+it+w,\tilde{\chi}) 
\Big(\frac{y^{w/2}-y^{w/4}}{w}\Big)^2 dw + O(1).
$$
Shifting contours to the left, this equals 
$$
-\sum_{\rho} \Big( \frac{y^{(\rho-1-it)/2}-y^{(\rho-1-it)/4}}{\rho-1-it}\Big)^2 
+O(1),
$$
where the sum is over all non-trivial zeros of $L(s,\tilde{\chi})$; the contribution of the 
trivial zeros may be absorbed into the $O(1)$ error term.  The contribution 
of zeros $\rho$ with $|\text{Im }\rho| >q$ is easily seen to be $\ll 
(\log q)/q \ll 1$. 
For a zero with $|\text{Im }\rho| \le q$ we see by our hypothesis that 
the numerator of our sum is $\ll y^{-j/(2\log q)}$.  Thus we obtain that 
$$
\sum_p \frac{\chi(p) \log p}{p^{1+it}} W(p) \ll y^{-j/(2\log q)} 
\sum_{|\text{Im }\rho|\le q} 
\frac{1}{|1+it-\rho|^2} 
+ 1. \tag{4.7}
$$

Now observe that if $\rho = \beta +i\gamma$ with $|\gamma|\le y$, 
and $\beta\le 1-j/\log q$ then
$$
\align
\frac{1}{|1+it-\rho|^2} &\ll \frac{1}{|1+1/\log q+it -\rho|^2} 
=\frac{1}{1+1/\log q -\beta} \text{Re }\frac{1}{1+1/\log q +it -\rho} 
\\
&\le \frac{\log q}{j} \text{Re }\frac{1}{1+1/\log q +it-\rho}.
\\
\endalign
$$
Thus the sum over zeros in (4.7)  is 
$$
\ll   \frac{\log q}{j} \sum_{\rho} \text{Re } \frac{1}{1+1/\log q+it-\rho} 
\ll \frac{(\log q)^2}{j},
$$
upon using a consequence of Hadamard factorization (see Davenport [4], chapter 12, equations (17) and (18)).  Since $y\ge q^{1/A}$ we deduce that 
$$
\sum_p \frac{\chi(p) \log p}{p^{1+it}} W(p) \ll e^{-j/(2A)} A^2 (\log y)^2 +1.
$$
Since $j\ge 4A \log A +D$ for a suitably large $D$, this proves the Lemma.

\enddemo

There are $\ll (\log q)^{10AC_2}$ characters $\chi \pmod q$ 
lying in $\Xi(j)$ for some $10A\log \log q \ge j\ge 4A \log A+D$.  
Therefore, by Proposition 4.1,
$$
\sum\Sb j= 4A\log A+D\endSb^{ 10A\log \log q} \sum\Sb \chi \in \Xi(j)\endSb
|\Psi(x,y;\chi,\Phi)| \ll \Psi(x,y;\chi_0,\Phi) \Big( 
(\log x)^{10AC_2+1} e^{-\sqrt{u}/20} + q^{-2}\Big).
$$
In our range for $x$ and $y$, we have $u\gg (\log \log x)^4$ 
and so we obtain that 
$$
\sum\Sb j= 4A\log A+D\endSb^{ 10A\log \log q} \sum\Sb \chi \in \Xi(j)\endSb
|\Psi(x,y;\chi,\Phi)| \ll \Psi(x,y;\chi_0,\Phi) \Big( \frac{1}{(\log x)^2} + \frac{1}{q}\Big). 
\tag{4.8}
$$


\head 5.  A few problem characters \endhead

\noindent In view of our work in \S 4,  it remains only to consider characters 
$\chi \pmod q$ with $\chi \in \Xi(j)$ for some $j\le 4A\log A +D$.  
By (2.12) there are only a bounded number $B=B(A)$, say, of such characters. 
We now define a set ${\Cal B}$ of {\sl problem characters}.  A non-principal character $\chi$ belongs 
to this set ${\Cal B}$ precisely if it has order at most $B$ and lies in 
$\Xi(j)$ for some $j\le 4A \log A +D$.   With Theorem 2 in mind, we define $H$ to be the 
subgroup of reduced residues $h\pmod q$ with $\chi(h)= 1$ for all $\chi \in {\Cal B}$.  
Since ${\Cal B}$ contains at most $B$ characters, all of order at most $B$, we see that 
the index of $H$ in $({\Bbb Z}/q{\Bbb Z})^*$ is at most $B^B$.

\proclaim{Proposition 5.1}  Retain our ranges for $x$, $y$, and $q$.  
If $\chi$ is non-principal, and $\chi \not\in {\Cal B}$  then 
$$
\Psi(x,y;\chi,\Phi) \ll \Psi(x,y;\chi_0,\Phi) \Big(\frac{1}{(\log x)^2 }+\frac{1}{q}\Big).
$$
\endproclaim 

\demo{Proof}  We first show that for all $|t|\le \sqrt{q}/(2B)$ 
$$
{\Bbb D}_{\alpha}(1,\chi(p)p^{-it};y)^2 \ge \frac{\sqrt{u}}{40B^2}. \tag{5.1}
$$
If not, there exists $t_{\chi}$ with $|t_{\chi}| \le \sqrt{q}/B$, 
and with 
$$
{\Bbb D}_\alpha(1,\chi(p)p^{-it_{\chi}};y)^2 \le 
\frac{\sqrt{u}}{40B^2}.  
$$
By the triangle inequality it follows that 
$$
{\Bbb D}_{\alpha}(1,\chi(p)^kp^{-ikt_{\chi}};y)^2 \le \frac{k^2 \sqrt{u}}{40B^2}.
$$
By Lemma 4.2 we see that $\chi$, $\chi^2$, $\ldots$, $\chi^{B+1}$ 
must all be in $\cup_{j\le 4A\log A+D} \Xi(j)$.  Since there are 
at most $B$ elements in $\cup_{j\le 4A\log A+D} \Xi(j)$, it follows that 
two of the $B+1$ characters listed above are the same.   But then 
$\chi$ would have order at most $B$ and would be in 
$\cup_{j\le 4A\log A+D} \Xi(j)$, contradicting our hypothesis that 
$\chi \not\in {\Cal B}$.

  We now use (2.11), invoking (5.1) for $|t|\le \sqrt{q}/(2B)$ 
and the rapid decay of ${\check \Phi}$ for $\sqrt{q}/(2B)\le |t|\le \sqrt{q}$.  We conclude 
that 
$$
\align
\Psi(x,y;\chi,\Phi) &\ll x^{\alpha} L(\alpha,\chi_0;y)\Big( \exp(-\sqrt{u}/(40B^2)) + q^{-2}\Big ) 
\\
&\ll \Psi(x,y;\chi_0,\Phi) \Big( (\log x) \exp(-\sqrt{u}/(40B^2)) + q^{-1}\Big).\\
\endalign
$$
Since $u \gg (\log \log x)^4$ in our ranges for $x$ and $y$, we obtain the Proposition. 
\enddemo

Now we examine more closely the situation for characters of bounded order.

\proclaim{Lemma 5.2}  Let $\chi$ be a character with order 
$k >1$.  In the range $|kt| \le 1/\log y$ we have 
$$
{\Bbb D}_{\alpha}(1,\chi(p)p^{-it};y) 
\gg \max\Big( {\Bbb D}_{\alpha}(1,\chi(p);y)^2,|t|^2 \log x\log y\Big),
$$
while in the range $1/\log y \le |kt| \le y$ we have 
$$
{\Bbb D}_{\alpha}(1,\chi(p)p^{-it};y) \gg \frac{u}{k^2 \log^2 u}.
$$
\endproclaim 
\demo{Proof} The triangle inequality gives
$$
{\Bbb D}_{\alpha}(1,\chi(p)p^{-it};y) \ge 
\frac{1}{k} {\Bbb D}_{\alpha}(1,p^{-ikt};y). \tag{5.2} 
$$
Consider first the case $y \ge |kt| \ge 1/\log y$.  Here we have
$$ 
{\Bbb D}_{\alpha}(1,p^{-ikt};y)^2 
\ge \frac{1}{\log y} \sum\Sb p\le y\\ p\nmid q \endSb \frac{1-\text{Re }p^{-ikt}}{p^{\alpha}} \log p 
\gg \frac{1}{\log y} \Big( \frac{y^{1-\alpha}}{1-\alpha} - \frac{y^{1-\alpha}}{|1-\alpha +ikt|}\Big),
$$
using the argument 
of the prime number theorem (using the Littlewood or 
Vinogradov zero-free regions for $\zeta(s)$).  Using (2.7) we 
conclude that 
$$
{\Bbb D}_{\alpha}(1,\chi(p)p^{-it};y)^2 \gg \frac{u}{k^2 \log^2 u}.  
$$
Using this in (5.2) we obtain our second assertion.

Now consider the range $|kt|\le 1/\log y$.  From (2.9) we obtain that 
$$
\sum\Sb p\le y\\ p\nmid q\endSb 
\frac{1-\text{Re }p^{-ikt}}{p^{\alpha}} 
\asymp \sum\Sb p\le y\\ p\nmid q\endSb 
\frac{(kt\log p)^2}{p^{\alpha}} \asymp
(kt)^2 \log x\log y. \tag{5.3}
$$
By (5.2) it follows that 
$$
{\Bbb D}_{\alpha}(1,\chi(p)p^{-it};y)^2 \gg t^2 \log x \log y,  
$$
which is one of the bounds in our first assertion.   
Moreover, by the triangle inequality we get that 
$$
{\Bbb D}_{\alpha}(1,\chi(p)p^{-it};y)+{\Bbb D}_{\alpha}(1,p^{it};y) \ge 
{\Bbb D}_{\alpha}(1,\chi(p);y).
$$
As in (5.3) we see that if $|t| \le 1/\log y$ then 
$$
{\Bbb D}_{\alpha}(1,p^{it};y)^2 \asymp t^2 \log x\log y.
$$
Therefore 
$$
{\Bbb D}_{\alpha}(1,\chi(p)p^{-it};y)^2 + O(t^2 \log x\log y) 
\gg {\Bbb D}_{\alpha}(1,\chi(p);y)^2,
$$
and the other bound claimed in our first assertion follows.

\enddemo


\proclaim{Lemma 5.3} Let $\chi$ be a character of order $1<k\le B$, 
and let $y\ge q^{\frac 1{4\sqrt{e}} +\delta}$.  Then 
$$
{\Bbb D}_{\alpha}(1,\chi;y)^2 \ge \frac{\delta}{4k} \log u + \log \delta +O(1).
$$
\endproclaim

\demo{Proof}  Put $z=y^{\sqrt{e}-\delta} > q^{\frac 14 +\delta}$.   We define the completely 
multiplicative function $f(n)$ by setting $f(p)=1$ for $p\le y$ and $f(p) =\chi(p)$ 
for $y<p\le z$.  Since $z<y^2$ we note that for $n\le z$ we have $f(n) 
= 1- \sum_{p|n} (1-f(p))$.  Therefore 
$$
\align
\text{Re }\sum_{n\le z} f(n) &= z - z \sum_{y\le p\le z}\text{Re } \frac{1-\chi(p)}{p} +o(z)
\ge z \Big(1 - \sum_{y\le p \le z}\frac 2p +o(1)\Big) 
\\
&= z\Big( 1- 2\log \frac{\log z}{\log y} +o(1)\Big),\\
\endalign
$$
and so 
$$
\Big| \sum_{n\le z} f(n)\Big| \gg \delta z. \tag{5.4} 
$$

Now let us write $f(n) = \sum_{d|n} g(d) \chi(n/d)$ 
where $g$ is a multiplicative function with $g(p)=1-\chi(p)$ 
for $p\le y$ and $g(p)=0$ for $y< p\le z$.  We see that 
$$
\sum_{n\le z} f(n) = \sum_{d\le z} g(d) \sum_{m\le z/d} \chi(m). \tag{5.5} 
$$
For the terms $d\le y^{\delta/2}$, so that $z/d>q^{\frac 14+\frac{\delta}{2}}$, 
we use Burgess's character sum estimates [3].  
The refinement of Heath-Brown (see [10], Lemma 2.4) applies, since our character 
has bounded order.   For such $d$ we see that $\sum_{m\le z/d} \chi(m) 
\ll (z/d)/(\log z)^3$ say, and hence 
$$
\sum_{d\le y^{\delta/2}} g(d) \sum_{m\le z/d} \chi(m) \ll 
\frac{z}{(\log z)^3} \sum_{d\le y^{\delta/2}} \frac{|g(d)|}{d} 
\ll \frac{z}{(\log z)^3} \sum_{d\le z} \frac{d(d)}{d} = o(z).
$$
The contribution of terms $d> y^{\delta/2}$ to (5.5) is 
bounded in magnitude by
$$
z\sum_{y^{\delta/2} \le d\le z} \frac{|g(d)|}{d} 
\le \frac{z}{y^{\delta(1-\alpha)/2}} \sum_{d\le z} \frac{|g(d)|}{d^{\alpha}} 
\ll \frac{z}{y^{\delta(1-\alpha)/2}} \exp\Big(\sum_{p\le y} \frac{|1-\chi(p)|}{p^{\alpha}}\Big).
$$
Since $\chi$ has order $k$, we have $|1-\chi(p)| \le k(1-\text{Re }\chi(p))$, and 
so the above is, using (2.7a,b), 
$$
\ll \frac{z}{u^{\delta/2}} \exp\Big(k {\Bbb D}_{\alpha}(1,\chi;y)^2\Big).
$$

From (5.4) and (5.5) we conclude that 
$$
\delta z \ll \Big| \sum_{n\le z} f(n) \Big| \ll o(z) + \frac{z}{u^{\delta/2}} \exp\Big( k{\Bbb D}_\alpha(1,\chi;y)^2\Big),
$$
and our Lemma follows.

\enddemo

\proclaim{Proposition 5.4}  Retain our ranges for $x$, $y$, and $q$.  
Let $\chi\pmod q$ be a character of order $1<k\le B$.  
Then, for any $1\le U \le \sqrt{u}$ we have, for some positive 
constant $c$,
$$
\align
\Psi(x,y;\chi,\Phi) &=\frac{1}{2\pi } \int_{|t|\le U/\sqrt{\log x\log y}} 
x^{\alpha+it} L(\alpha+it,\chi;y) {\check \Phi}(\alpha+it)dt \\
&\hskip .5 in+ O\Big( \Psi(x,y;\chi_0,\Phi)\Big(\frac{1}{q^2}+
\frac{1}{(\log x)^2} + e^{-cU^2}\Big)\Big) .
\\
\endalign
$$
If $A<4\sqrt{e}-100\delta$ then for some small positive constant $c$ 
$$
\Psi(x,y;\chi,\Phi) \ll \Psi(x,y;\chi_0,\Phi) u^{-c\delta}.
$$
\endproclaim 

\demo{Proof}  We start with the expression (2.11).  We split the integral over 
$|\text{Im }s|\le \sqrt{q}$ into various ranges.  The rapid decay of ${\check \Phi}(s)$ 
shows that the contribution to the integral from $y\le |\text{Im }s|\le \sqrt{q}$ is 
$\ll \Psi(x,y;\chi_0,\Phi)q^{-2}$.  In the range $1/\log y \le |\text{Im }s|\le y$ 
we use the second bound of Lemma 5.2.  Thus the contribution of this range 
is $\ll \Psi(x,y;\chi,\Phi) (\log x) \exp(-Cu/\log^2 u)$ for 
some constant $C$.  Since $u\gg (\log \log x)^4$ this contribution is $\ll 
\Psi(x,y;\chi,\Phi)/(\log x)^2$.  In the range $U/\sqrt{\log x\log y} 
\le |\text{Im }s| \le 1/\log y$, we use the first assertion of Lemma 5.2 which 
gives ${\Bbb D}_{\alpha}(1,\chi(p)p^{-it};y)^2 
\gg |t|^2 \log x\log y$.  It follows that the contribution of this range is 
$\ll x^{\alpha} L(\alpha,\chi_0;y) \exp(-cU^2)/\sqrt{\log x\log y} 
\ll \Psi(x,y;\chi_0,\Phi) e^{-cU^2}$  for some positive constant $c$.  Piecing 
these statements together, we obtain the first assertion of the Proposition.

To prove the second assertion, we choose $U=\sqrt{c^{-1} \log u}$.  From Lemmas 
5.2 and 5.3 it follows that for $|t| \le U/\sqrt{\log x \log y}$
$$
{\Bbb D}_{\alpha} (1,\chi(p)p^{-it};y)^2 \gg {\Bbb D}_{\alpha}(1,\chi(p);y)^2 
\gg \delta \log u.
$$
Using this estimate to handle the integral in our first assertion, the Proposition follows. 

\enddemo

\head 6.  Proofs of the main Theorems \endhead

\demo{Proof of Theorem 1}  Combining (2.1), (3.2), (4.8), Proposition 5.1 and the second part 
of Proposition 5.4 we obtain that if $A \le 4\sqrt{e}-100\delta$ then 
$$
\Psi(x,y;q,a,\Phi) = \frac{1}{\phi(q)} \Psi(x,y;\chi_0,\Phi) \Big(1 + O(u^{-c\delta})\Big),
$$
for some positive constant $c$.  We now take $\Phi$ to be $1$ on $[0,1-\epsilon]$ 
and $0$ on $[1,\infty)$ to get a lower bound for $\Psi(x,y;q,a)$; and 
$\Phi$ to be $1$ on $[0,1]$ and $0$ on $[1+\epsilon,\infty)$ to get an upper 
bound for $\Psi(x,y;q,a)$.  Taking $\epsilon$ sufficiently small, we obtain Theorem 1. 

\enddemo 

\demo{Proof of Theorem 2} Combining (2.1), (3.2), (4.8) and Proposition 5.1 we obtain that 
$$
\Psi(x,y;q,a,\Phi) = 
\frac{1}{\phi(q)} \Psi(x,y;\chi_0,\Phi)\Big( 1+O\Big(\frac{1}{q}+\frac{1}{(\log x)^2}\Big)\Big) +
\frac{1}{\phi(q)}\sum\Sb \chi \in {\Cal B}\endSb \overline{\chi(a)} \Psi(x,y;\chi,\Phi).\tag{6.1}
$$
Once again we take $\Phi$ to be $1$ on $[0,1-\epsilon]$ 
and $0$ on $[1,\infty)$ to get a lower bound for $\Psi(x,y;q,a)$; and 
$\Phi$ to be $1$ on $[0,1]$ and $0$ on $[1+\epsilon,\infty)$ to get an upper 
bound for $\Psi(x,y;q,a)$.  Note that in either case ${\check \Phi}(\alpha+it) 
=1/(\alpha+it) + O(\epsilon)$.   Therefore, from (6.1), (2.9) and Proposition 5.4
(taking there $U=1/\sqrt{\epsilon}$) we may conclude that 
(for large $x$, $y$ and $q$ lying in our ranges) 
$$
\align
\Psi(x,y;q,a) &= \frac{1}{\phi(q)} \frac{x^{\alpha} L(\alpha,\chi_0;y)}{\sqrt{2\pi \phi_2(\alpha,\chi_0;y)}} 
+O\Big(\frac{\sqrt{\epsilon}}{\phi(q)} \Psi(x,y;\chi_0)\Big)\\
&+ \frac{1}{\phi(q)} \sum\Sb \chi \in {\Cal B}\endSb \overline{\chi(a)} 
\frac{1}{2\pi} \int_{|t|\le 1/\sqrt{\epsilon 
\log x\log y}} x^{\alpha+it} L(\alpha+it,\chi;y) \frac{dt}{\alpha+it} 
.\tag{6.2}\\
\endalign
$$

Recall that $H$ is the subgroup of residues $h$ such 
that $\chi(h) = 1$ for all $\chi\in {\Cal B}$, and that it has index 
at most $B^B$.   
If $a/b\in H$ then $\chi(a)=\chi(b)$ for each $\chi \in {\Cal B}$.  
Therefore (6.2) gives identical expressions for both $\Psi(x,y;q,a)$ and 
$\Psi(x,y;q,b)$.  Consequently 
$$
\Psi(x,y;q,a) = \Psi(x,y;q,b) + O\Big(\frac{\sqrt{\epsilon}}{\phi(q)} \Psi(x,y;\chi_0)\Big).
$$
This proves Theorem 2. 

\enddemo

\enddocument
\head 3.  My version of Halasz's method \endhead

\noindent Let $f$ be a multiplicative function with $|f(n)|\le 1$ 
for all $n$ and put $F(s)=\sum_{n=1}^{\infty} f(n) n^{-s}$.  We seek to 
estimate 
$$
\sum_{n} f(n) \log n \Phi(n/x),
$$
where $\Phi$ is a nice smooth function supported in $[1,2]$.  The 
trivial bound is $\ll x\log x$, and we want an improvement if 
$f(p)$ doesn't look too much like $1$ or $p^{i \beta}$ for some 
$\beta$. Note that by Mellin inversion we have 
$$
\sum_{n} f(n) \log n \Phi(n/x) 
= -\frac{1}{2\pi i} \int_{c-i\infty}^{c+i\infty} F^{\prime}(s) 
x^s {\check \Phi}(s) ds.
$$

\noindent Recall that there is a constant $C>0$ such that 
if $1 \le \psi \le \log q$ then the region 
$$
{\Cal R}_{\psi} := \{ \sigma > 1-\psi/\log q , \ \ |t| \le q \}
$$
contains fewer than $e^{C\psi}$ zeros of 
$\prod_{\chi\pmod q} L(s,\chi)$.  Thus for at most $e^{C\psi}$ 
characters can there be a zero in ${\Cal R}_{\psi}$.

\proclaim{Lemma 2.1} Suppose that  
the region ${\Cal R}_{\psi}$ is free of zeros of $L(s,\chi)$.   
Then, we have, for $|t| \le \log q$ (say) and $s=c+it$
$$
\log |L(s,\chi;y)| \ll .
$$
\endproclaim
\demo{Proof} We begin with a bound for 
$$
\sum_{p\le z} \frac{\chi(p)\log p}{p^{it}}.
$$
Using the explicit formula we see that this is 
$$
= -\sum_{|\rho|\le q} \frac{z^{\rho}}{\rho} + O\Big( \frac{z\log^2 z}{q}\Big) 
\ll  z^{1-\psi/\log q} \log^2 q+ \frac{z\log^2 z}{q}.
$$
Trivially the sum is also $\ll z$.

\enddemo

\enddocument